# Saccheri`s Rectilinear Quadrilaterals

Prodromos Filippidis
Glyka Nera Attiki GR-15354, Greece
E- mail: prodromos.filippidis@gmail.com

*Abstract*— We study Saccheri`s three hypotheses on a two right-angled isosceles quadrilateral, with a rectilinear summit side. We claim that in the Hilbert`s foundation of geometry the euclidean parallelism is a theorem, and in the h-plane the hyperbolic parallelism under a hyperbolic transformation has image the euclidean parallelism. We prove that Saccheri`s rectilinear quadrilaterals can be only rectangle. Finally we believe that the independence of the euclidean parallel postulate is just a matter of philosophy of logic.

*Keywords*— rectilinear quadrilaterals. Euclidean parallelism. Hyperbolic parallelism. Rectilinear bisector

## I. Introduction

We study Saccheri`s three hypotheses on a two right-angled isosceles quadrilateral, with a rectilinear summit side. We claim that in the Hilbert`s foundation of geometry the euclidean parallelism is a theorem. Also we examine the inversion of the hyperbolic parallel lines, in the h-plane (the upper half euclidean plane). The paper is organized in the following way, at section 2 we present, the used definitions, postulates, and theorems. At section 3 we prove that, if for Saccheri`s rectilinear quadrilaterals the hypothesis of the acute angle, or the hypothesis of obtuse angle is true. Then a rectilinear isosceles triangle there exists, in which the perpendicular bisector to the congruent sides and the rectilinear bisector of the vertex angle are not intersecting. At section 4 we prove that the rectilinear bisector of the vertex angle of a rectilinear isosceles triangle and the perpendicular bisector of the congruent sides of this triangle are intersecting.

## II. Used definitions, postulates, and theorems

In this paper we consider a Hilbert plane, satisfying the axioms (I1)-(I3), (B1)-(B4), (C1)-(C6), the axiom of continuity, with the definitions and theorems deduced by them. [1]. Also: To prove the theorem of exterior angle we double the median BD of any triangle ABC by drawing the circle (D, DB). To prove the theorem that the rectilinear perpendicular to a straight line through a point not on the line is unique, we use the theorem of exterior angle and the Hilbert`s theorem that all right angles are congruent to one another.

## III. Theorem I

If for Saccheri`s rectilinear quadrilaterals the hypothesis of the acute angle, or the hypothesis of the obtuse angle is true. Then a rectilinear isosceles triangle there exists, in which the perpendicular bisector to the congruent sides and the rectilinear bisector of the vertex angle are not intersecting. Proof: For the Saccheri`s quadrilateral ΑΕΖΔ and its common perpendicular ΓΜ to the base ΔΖ and to the summit ΑΕ, is valid (a) ΑΔ=ΕΖ>ΓΜ if the hypothesis of the acute angle is true "Fig. 1" and (b) ΑΔ=ΕΖ<ΓΜ if the hypothesis of the obtuse angle is true "Fig. 2". Also, on the straight line ΜΓ there exists a point Β, so that ΑΔ=ΕΖ=ΒΓ, to the right of the point Μ or to the left of the point Μ, if the case (a) or the case (b) is respectively valid. The point Β together with the side ΑΕ of the quadrilateral ΑΕΖΔ forms the rectilinear isosceles triangle ΑΒΕ. Since the congruent sides ΑΒ and ΒΕ of the triangle ΑΒΕ are also the summits of the acute-angled "Fig. 1" or obtuse-angled "Fig. 2" isosceles quadrilaterals ΑΒΓΔ and ΒΕΖΓ it implies that the perpendicular bisectors ΠΟ and ΤΣ to the congruent sides of the isosceles triangle ΑΒΕ, as well the bisector ΒΓ of the vertex angle ΑΒΕ are perpendiculars to the base ΔΖ of the acute-angled or the obtuse-angled quadrilateral ΑΕΖΔ, and since the perpendicular to a line through a point not on the line is unique, we can conclude that they are not intersecting. In the section 4 we shall use the theory of function to examine if the above conclusion, that the rectilinear perpendicular bisectors to the congruent sides of a rectilinear isosceles triangle are not intersecting, is true.

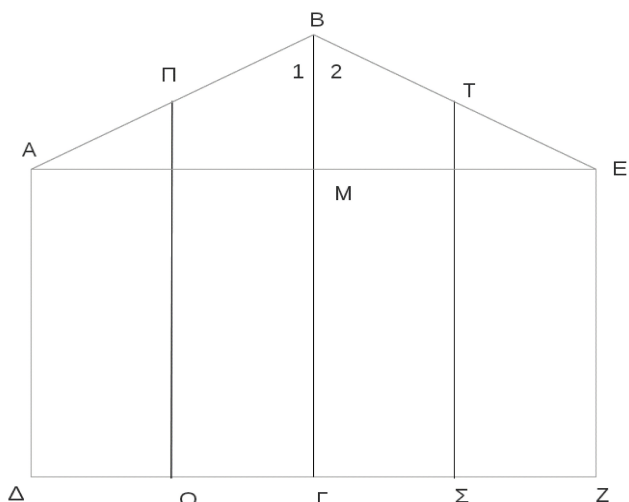

Fig. 1

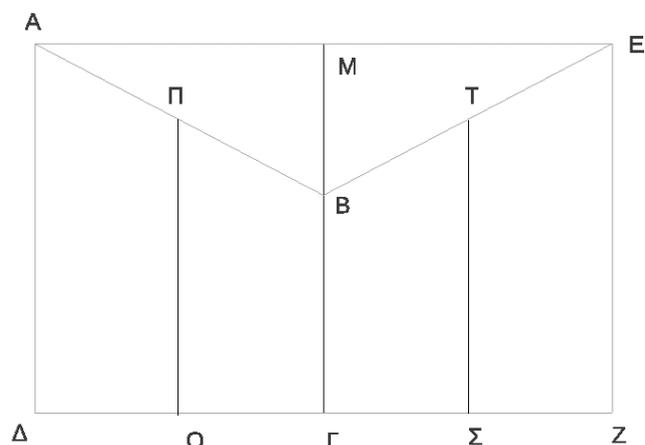

Fig. 2

IV. THEOREM II

The rectilinear perpendicular bisectors to the congruent sides of a rectilinear isosceles triangle are intersecting [2]

We consider the Cartesian system "Fig. 3" where the x-axis coincides with the base AE of the isosceles triangle ABE of figure 1 and the y-axis coincides with the bisector BΓ of the vertex angle ABE of it. The points Π(a1, b1) and Π2(a2, 0) absolutely determine the straight line ΠΟ which is perpendicular bisector to the side AB of the isosceles triangle ABE, and the points T(a1', b1) and T2(a2', 0) determine the straight line TΣ which is perpendicular bisector to the side BE of the same isosceles triangle ABE. Since straight lines absolutely determined by their two points (incidence axioms I1-I3) their functions are (1-1) and for, if they are not perpendiculars to the x-axis or to the y-axis. The points Π and T have the same ordinate (y-value) b1 since the straight line ΠT is perpendicular to the y-axis, because it is summit of the isosceles Saccheri rectilinear quadrilateral ΠTΣΟ The perpendiculars to the base ΔZ straight lines ΠΟ and TΣ are congruent because they are corresponding parts of the congruent triangles ΔΠΟ and ZTΣ. Now we shall compare (in three cases) the functions of the straight lines ΠΟ and TΣ using the points which determine them.

Case1: A function f is equal to a function g when all have the same (set of) definition, the same set values and assign equal arguments to equal values:

$$f(a)=b \text{ iff } g(a)=b$$

From points Π, Π2 and T, T2 we have that f(a1)=b1 and g(a1')=b1 for each $a1 \neq a1'$ and b1=b1 This implies that the functions f and g are not equal.

Case 2: The (foreign) compound of two functions f: A → B and g: A' → B', where A, A' are disjoints sets, is the mapping $f \cup g$ : $A \cup A'$ → $B \cup B'$ who define as $f \cup g(a)$ =f(a) and $f \cup g(a')$ =g(a') for each $a \in A$, $a' \in A'$

From points Π, Π2 and T, T2 we have $f \cup g(a1)$ =f(a1)=b1, and $f \cup g(a1')$ =g(a1')=b1 for $a1 \neq a1'$ This also implies that the functions f and g are not foreign. Since we reject the cases 1 and 2 we have to accept the case 3. That is, that the functions f and g are intersecting at a point.

Case 3: The intersection of two functions f: A → B and g: A' → B' is the mapping $f \cap g$ : $A \cap A'$ → $B \cap B'$ defined as $f \cap g(a)$ =b iff f(a)=g(a)=b for each $a \in A \cap A'$

It is easy proved that the point of intersection of the perpendicular bisectors to the congruent sides of the isosceles triangle ABE lies on the y-axis, which coincides with the bisector BΓ, of the vertex angle ABE, and its (x-value) is 0. This means that the point of intersection belongs to the set of arguments (x-values) of functions f and g that is: $f \cap g(0)$ =b since f(0)=g(0)=b for $0 \in A \cap A'$

The acceptance of the case 3 means that the perpendicular bisectors to the sides of an isosceles triangle are intersecting. This last implies that the hypotheses that the summit angles of a rectilinear Saccheri quadrilateral are acute or obtuse angles are not true.

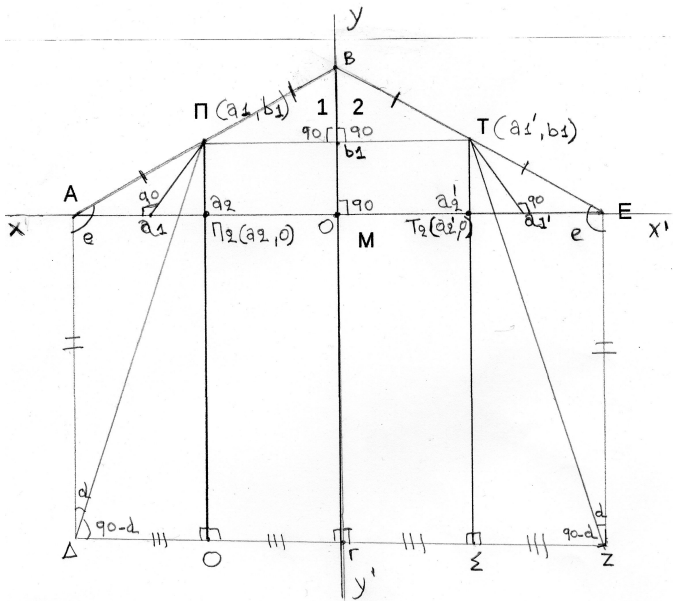

Fig. 3

## V. CONCLUSIONS

Since we reject the hypotheses of the acute angles and the obtuse angles, the acceptance that the rectilinear Saccheri`s quadrilaterals are only rectangle quadrilaterals, implies that the euclidean parallelism is proved as a theorem in the Hilbert foundation of geometry, and as that it can be used. That is, all the euclidean straight lines, which are in h-plane, the upper half euclidean plane, and they are parallel in the euclidean sense are not met at a point at infinity. This last result, that the hyperbolic motion inversion, which takes place, in the h-plane (the upper half euclidean plane) of parallel hyperbolic straight lines (euclidean half circles with center on u-line) transforms them into euclidean lines, perpendicular to u-line which are euclidean parallel straight lines, since they are not met at a point at infinity. That is, in the h-plane the hyperbolic parallelism under a hyperbolic transformation has image the euclidean parallelism. We think that the proof of euclidean parallel postulate as theorem has not any other result, except the above, on the parallelism, on the hyperbolic geometry. We believe that the independence of the euclidean parallel postulate is just a matter of philosophy of logic. The proof that Saccheri rectilinear quadrilaterals can be only rectangle implies that from these quadrilaterals are not deduced any conclusion for non-euclidean geometries. An acute-angled Saccheri quadrilateral of hyperbolic geometry, has only curved summit side.